\documentclass[fleqn, eqnumis, pdf]{hjms}
\usepackage{amsmath,amsthm,amscd,amsfonts,amssymb,enumerate,tikz-cd,mathrsfs,array}
\usepackage{graphicx}		
\usepackage{color}
\usepackage[colorlinks]{hyperref}
\usepackage[T1]{fontenc}
\usepackage{footnote}
\usepackage{multicol}
\usepackage{booktabs}
\usepackage[flushleft]{threeparttable}
\usepackage[font=small,labelfont=bf]{caption}

\newcommand{\bte}{\begin{theorem}\quad  }
\newcommand{\ete}{\end{theorem} }
\newcommand{\bpr}{\begin{proposition}\quad  }
\newcommand{\epr}{\end{proposition} }
\newcommand{\ble}{\begin{lemma}\quad }
\newcommand{\ele}{\end{lemma}}
\newcommand{\bco}{\begin{corollary}\quad }
\newcommand{\eco}{\end{corollary} }
\newcommand{\bex}{\begin{example}\quad \rm }
\newcommand{\eex}{\end{example} }
\newcommand{\bde}{\begin{definition}\quad \rm }
\newcommand{\ede}{\end{definition} }
\newcommand{\brm}{\begin{remark} \quad \rm}
\newcommand{\erm}{\end{remark} }
\newcommand{\bpf}{\begin{proof}[\bf{Proof.\quad}] \rm}
\newcommand{\epf}{ \end{proof}}
\newcommand{\bdm}{\begin{displaymath} }
\newcommand{\edm}{\end{displaymath} }
\newcommand{\be}{\begin{eqnarray*}}
\newcommand{\ee}{\end{eqnarray*}  }

\newcommand{\lb}{\label}

\newcommand{\ba}{\begin{align*}}
\newcommand{\ea}{\end{align*}}

\makeatletter
\newcommand\cupop{\mathop{\operator@font \cup}\nolimits}
\makeatother

\begin{document}

\title{Monoids over which products of indecomposable acts are indecomposable}

\author{ Mojtaba Sedaghatjoo \footnote{ Department of Mathematics, College of Sciences,
 Persian Gulf University, Bushehr, Iran,  \ Email: {\tt sedaghat@pgu.ac.ir }}\footnote{Corresponding Author.} and
Ahmad Khaksari\footnote{ Department of Mathematics, Payame Noor University, Tehran, Iran,  \ Email: {\tt a\_khaksari@pnu.ac.ir }}}
{ }

\begin{abstract}
In this paper we prove that for a monoid $S$, products of indecomposable right $S$-acts are indecomposable if and only if $S$ contains a right zero. Besides, we prove that subacts of indecomposable right $S$-acts are indecomposable if and only if $S$ is left reversible. Ultimately, we prove that the one element right $S$-act $\Theta_S$ is product flat if and only if $S$ contains a left zero.
\end{abstract}

\begin{keyword}
Indecomposable act, left reversible monoid,  Baer criterion, product flat, super flat.
\end{keyword}

\begin{AMS}
Primary: 20M30; Secondary: 20M50.
\end{AMS}

\section{\bf Introduction}

Throughout this paper, $S$ stands for a monoid and $1$ denotes its identity element. A nonempty set $A$ together with a mapping $A \times S \rightarrow A,~(a,s)\rightsquigarrow as$, is called a right $S$-act or simply an act (and is denoted by $A_S$) if $a(st)=(as)t$ and $a1=a$ for all $a\in A,~s,t \in S$. Left $S$-acts can be defined similarly. We mean by $A\sqcup B$ the disjoint union of  sets $A$ and $B$. The one element act is called zero act and is denoted by $\Theta_S$. A right $S$-act $A_S$ is called decomposable provided that there exist subacts $B_S,C_S\subseteq A_S$ such that $A_S=B_S \cup C_S$ and $B_S \cap C_S= \emptyset$; in this case $A_S=B_S \cup C_S$ is called a decomposition of $A_S$. Otherwise $A_S$ is called indecomposable. For a nonempty set $I$, $S^I$ denotes the set $\prod \limits _I S$, endowed with the natural componentwise right $S$-action ${(s_i)}_{i\in I}s={(s_is)}_{i\in I}$. We refer the reader to \cite{adamek,Kilp} for more details on the concepts mentioned in this paper.

 Since for a given monoid $S$ any right $S$-act $A_S$ is uniquely the disjoint union of indecomposable acts called  indecomposable components of $A_S$, analogous to the bricks forming a wall, indecomposable acts deserve to be taken into consideration. A pioneering work in this account goes back to \cite{niko}, where the collection of all indecomposable right $S$-acts is partitioned into equivalence classes corresponding
to the components of the right $S$-act $\mathscr{R}$ formed by letting $S$ act on its right congruences by translation.

 As mentioned, every right $S$-act $A_S$ has a unique decomposition into indecomposable subacts, indeed, indecomposable components of $A_S$ are the equivalence classes of the relation $\sim $ on $A_S$  defined in \cite{Ren} by $a\sim b$ if there exist $s_1,s_2, \ldots ,s_n,t_1,t_2, \ldots ,t_n \in S,~a_1,a_2, \ldots , a_n \in A_S$ such that
\[ a=a_1s_1,~ a_1t_1=a_2s_2,~ a_2t_2=a_3s_3, \ldots,a_{n-1}t_{n-1}=a_ns_n,~a_nt_n=b. \] We shall call this sequence of equalities a scheme of length $n$.  Therefore, elements $a,b\in A_S$ are in the same indecomposable component if and only if there exists a scheme of length $n$  as above connecting $a$ to $b$. Note that for a natural number $m>n$, the scheme length can be increased to  $m$ by adding the equality $b.1=b.1$ to the end of scheme iteratively.

The paper comprises three sections as follows. In the first section we present a short account of the needed notions. The second one concerns with indecomposable acts over left reversible monoids. We prove that in Baer criterion for acts, the condition of possessing a zero element can be abandoned in case that $S$ is not left reversible. In third section we engage in the main results of this paper, that is, conditions under which
properties of indecomposability, product flatness and super flatness are preserved under products. Furthermore, we prove that for the one element act $\Theta_S$, the tensoring functor  $\Theta_S \otimes -$ preserves limits if and only if it preserves products, equivalently, products of indecomposable left $S$-acts are indecomposable.

\section{\bf Indecomposable acts over left reversible monoids}

In this section we investigate indecomposable acts over left reversible monoids (that are monoids satisfying nonempty intersection for any pair of right ideals) and give some characterizations for left reversible monoids regarding indecomposability property. In the next proposition we show that for left reversible monoids the length of the preceding scheme can be considered to be $2$.

\bpr \label{pr1}For a monoid $S$ the following are equivalent:

$~~~~\,i)$ $S$ is a left reversible monoid,

$~~~\,ii)$ a right $S$-act $A_S$ is indecomposable if and only if for any $a,a'\in A_S$ there exist $s,s' \in S$ such that $as=a's'$,

$~~\,iii)$ any indecomposable right $S$-act contains at most one zero element.\epr

\bpf $i\Longrightarrow ii$. Let $S$ be a left reversible monoid and suppose that $a,a'\in A_S$ for an indecomposable right $S$-act $A_S$. So there exists a scheme connecting $a$ to $a'$, of the form \begin{center} $ a=a_1s_1,~ a_1t_1=a_2s_2,~ a_2t_2=a_3s_3, \ldots,~a_nt_n=a',$ \end{center}  for $a_i\in A_S,~s_i,t_i\in S,~1\leqslant i \leqslant n $. Left reversibility of $S$ provides  $u_1,v_1 \in S$ such that $s_1u_1=t_1v_1,$ and in consequence $au_1=a_1s_1u_1=a_1t_1v_1=a_2s_2v_1$. Proceeding inductively, we get $u,v\in S$ satisfying $au=a_nt_nv=a'v$ as desired.

$ii\Longrightarrow iii$. On the contrary suppose that an indecomposable right $S$-act contains two different zero elements namely $\theta_1$ and $\theta_2$. By assumption there exist $s,t\in S$ for which $\theta_1=\theta_1s=\theta_2t=\theta_2$ a contradiction.

$iii\Longrightarrow i$. By way of contradiction suppose that $I\cap J=\emptyset$ for two right ideals $I$ and $J$ of $S$. Now define a right congruence $\rho$ on $S$ by $x\rho y$ if $x=y$ or $x,y\in I$ or $x,y\in J$. Take $a\in I$ and $b\in J$. Then $S/\rho$ is a cyclic indecomposable right $S$-act with two different zero elements namely $[a]$ and $[b]$, a contradiction.
\epf

Recall that Baer criterion for right $S$-acts asserts that a right $S$-act is injective if and only if it possesses a zero element and is injective relative to all inclusions into cyclic right $S$-acts. In what follows we prove that if $S$ is not left reversible then the condition of possessing a zero element in Baer criterion could be omitted.
\bpr \lb{pr8} Let $S$ be a monoid that is not left reversible. A right $S$-act $Q_S$ is injective if and only if it is injective relative to all inclusions into cyclic right $S$-acts.
\epr
\bpf Necessity is clear. To prove sufficiency, let $Q_S$ be a right $S$-act that is injective relative to all inclusions into cyclic right $S$-acts. Suppose that $I\cap J=\emptyset$ for two right ideals $I$ and $J$ and $\rho$ is the Rees congruence on $S_S$ defined by the right ideal $J$. Consider the homomorphism $f:I_S \longrightarrow Q_S$ given by $f(i)=qi,i\in I$ for some $q\in Q_S$. Since $I$ can be identified with a subact of $S/\rho$, our assumption yields a homomorphism $\bar{f}: S/\rho \longrightarrow Q_S$ making the following diagram commutative.
 \begin{center}
\begin{tikzcd}
I_S \arrow{r}{\subseteq} \arrow{d}[swap]{f}
&S/\rho \arrow{ld}{\bar{f}} \\
Q_S &
\end{tikzcd}
\end{center}
Now $S/\rho$ contains a zero element and so does $Q_S$ which thanks to Baer criterion $Q_S$ is injective.
\epf

Here a question can be posed that \\
{\bf whether a monoid $S$ over which injective acts are precisely ones that are injective relative to all inclusions into cyclic acts, is not left reversible.}

In the next proposition we characterize monoids over which subacts of indecomposable acts are indecomposable.

\bpr \lb{pr2}For a monoid $S$ all subacts of indecomposable right $S$-acts are indecomposable if and only if $S$ is left reversible. \epr

\bpf \textit{Necessity}. Let $a,b\in S$. Since $S$ is indecomposable, our assumption implies that $aS\cup bS$ is indecomposable and therefore $aS \cap bS \neq \emptyset$.

\textit{Sufficiency}. This is a straightforward result of Proposition \ref{pr1}, part $ii)$.
\epf

Recall that for a nonempty set $I$, $I^S$ is an $|I|$-cofree right $S$-act where $fs$ for $f\in I^S, ~s\in S$ is defined by $(fs)(t)=f(st)$ for every $t\in S$. It should be mentioned that the 1-cofree objects or terminal objects in {\bf Act-}$S$ are precisely one element acts which are indecomposable.  The next proposition characterizes monoids over which non-zero cofree acts are decomposable.
\bpr \lb{pr3} For a monoid $S$ the following are equivalent:

$~~~~\,i)$ all non-zero cofree $S$-acts are decomposable,

$~~~\,ii)$ there exists a non-zero decomposable cofree right $S$-act,

$~~\,iii)$ $S$ is left reversible.

\epr
\bpf $i\Longrightarrow ii$ is clear. $ii\Longrightarrow iii$. By way of contradiction suppose that $aS \cap bS=\emptyset$ for some $a,b\in S$. Let $X^S$ be a non-zero decomposable $|X|$-cofree act and $f,g \in X^S$.  Let $h \in X^S$ be given by $$h(x)=  \begin{cases} f(x) & \text{if } x\in aS, \\ g(x) &  \text{otherwise.}
\end{cases}$$
So we get the scheme $f=f.1,~fa=ha,~hb=gb,~g.1=g$, which implies that $f$ and $g$ are in the same indecomposable component. Therefore $X^S$ is indecomposable a contradiction.

$iii\Longrightarrow i$. Let $S$ be a left reversible monoid and $X^S$ be a non-zero cofree $S$-act. Take constant functions $f=c_{x_1}$ and $g=c_{x_2}$ in $X^S$ for different elements $x_1$ and $x_2$ in $X$. Then $f$ and $g$ are zero elements of $X^S$ (note that zero elements of $X^S$ are precisely constant functions). In light of Proposition \ref{pr1}, since $X^S$ contains two different zero elements, it is not indecomposable.
\epf

Regarding the fact that any scheme in an arbitrary act  can be translated into another one in its factor act, generally factor acts of indecomposable acts are indecomposable. It is clear that coproducts of indecomposable acts are not indecomposable, though the next proposition states that pushouts of indecomposable acts are indecomposable.

It is worth pointing out that for a monoid $S$, since right $S$-acts are nonempty, the category of right $S$-acts is not complete nor cocomplete. Indeed, this category has products and coequalizers and has neither coproducts and equalizers. Note that in this category coproduct of nonempty families of objects exists. Hence, coproducts can not be considered as a sort of pushouts.

\bpr \lb{pr7} For a monoid $S$ consider a pushout situation
\vspace{.5cm}
\begin{center} $\begin{CD}
Y_1 \\
@A{f_1}AA \\
X_S @>>{f_2}> Y_2
\end{CD}$
\end{center}
\vspace{.5cm}
in the category {\bf Act}-$S$ of right $S$-acts where $Y_1$ and $Y_2$ are indecomposable and suppose that $((q_1,q_2),Q_S)$ is the pushout of the pair $(f_1,f_2)$. Then $Q_S$ is indecomposable. \epr
\bpf It is known that $Q_S$ is isomorphic to $(Y_1\sqcup Y_2)/\nu$ where $\nu$ is the congruence relation on $(Y_1\sqcup Y_2)$ generated by all pairs $(f_1(x),f_2(x)),x\in X$. Let $[y_1],[y_2]\in (Y_1\sqcup Y_2)/\nu$ for some $y_1,y_2\in  Y_1\sqcup Y_2$. Since $Y_1$ and $Y_2$ are indecomposable, in view of the preceding argument we should just engage in the case that $y_1\in Y_1,y_2\in Y_2$ or vice versa. Without restriction of generality, we can consider only the first case. Take an element $x_0\in X$. Therefore there exist two schemes connecting $y_1$ to $f_1(x_0)$ and  $y_2$ to $f_2(x_0)$. Thus we get two schemes in $(Y_1\sqcup Y_2)/\nu$ connecting $[y_1]$ to $[f_1(x_0)]$ and $[y_2]$ to $[f_2(x_0)]$ which, using the equality $[f_1(x_0)]=[f_2(x_0)]$ provides the desired result.
\epf
As amalgamated coproduct of objects in a category is a sort of pushout, the next corollary follows.
\bco \lb{co3} Amalgamated coproducts of indecomposable acts are indecomposable.
\eco
\section{\bf Products of indecomposable acts}

There have been published several works on preservation of acts properties under products, for instance \cite{Bul1,Bul,Bul2,Gould,sedaghat}. In this section we investigate another version of the problem for indecomposability property. Note that products of indecomposable acts are not indecomposable in general, for instance if $S$ is a left zero semigroup with an identity element externally adjoined, then there is no scheme in $S\times S$ connecting $(1,a)$ to $(a,1)$ for $ 1\neq a\in S$. As a product of a family of right $S$-acts is a sort of pullback, then indecomposability property is not preserved under pullback and consequently coamalgamated product. It is easy to check that for non-empty sets $I$ and $J$ with $|J|\leq |I|$, if $S^I$ is indecomposable, then so is $S^J$. Now suppose that $S^{S\times S}$ is indecomposable. Let $I$ be a non-empty set and  $(a_i)_I,(b_i)_I \in S^I$. Put $J=\{(a_i,b_i)|~i\in I\}$. We index $J$ by a set $K$ as $J=\{(u_k,v_k)|~k\in K\}$. Since $S^{S\times S}$ is indecomposable and $|K|=|J| \leq |S\times S|$, by the preceding argument, $S^K$ is indecomposable and hence there exists a scheme in $S^J$ connecting $(u_k)_K$ to $(v_k)_K$. The corresponding scheme in $S^I$ is the one connecting $(a_i)_I$ to $(b_i)_I$ as desired. Thereby, the next corollary is obtained.

 \bco \lb{co4} For a  monoid $S$, $S^{S\times S}$  is indecomposable if and only if $S^I$ is indecomposable for each nonempty set $I$.
 \eco
 A subject of interest in the study of tensor products is preservation of limits by tensoring functor $A_S\otimes -$ for a right $S$-act $A_S$ which is investigated  in \cite{Bul}. Following terms used  in this reference, a right $S$-act $A_S$ is called (finitely) super flat if the functor $A_S\otimes -$ preserves all (finite) limits, and (finitely) product flat if it preserves all (finite) products. Now if finite products of indecomposable acts are indecomposable then $S\times S$ is indecomposable. In the next theorem we show that this is a sufficient condition for finite products of indecomposable acts to be indecomposable which is equivalent to the condition that the one-element left $S$-act $_S\Theta$ is finitely product flat. Besides in the sequel we show that products of indecomposable acts are indecomposable if and only if the one element left $S$-act $_S\Theta$ is product flat.
\bte \lb{th1} For a monoid $S$ the following are equivalent.

 $~~i)$ finite products of indecomposable acts are indecomposable,

 $~ii)$ finite products of cyclic acts are indecomposable,

 $iii)$ ${S^n}$ is indecomposable for each $n \in \mathbb{N}$,

 $\,iv)$ ${S^n}$ is indecomposable for some $1\neq n \in \mathbb{N}$,

 $\,~v)$ $S \times S$ is indecomposable,

 $\,vi)$ the one element left $S$-act $_S\Theta$ is finitely product flat.
\ete
 \bpf It is sufficient to prove the implication $v\Longrightarrow i$ and the term $v \Longleftrightarrow vi$ is valid by \cite[Corollary 2.5]{Bul}.  Suppose that $S\times S$ is indecomposable. We just need to prove that the product of two indecomposable right $S$-act is indecomposable and then applying induction provides the desired result. Let $A_S$ and $B_S$ be indecomposable right $S$-acts and let $(a,b),(a',b')\in A \times B$ for some $a,a'\in A_S,~b,b'\in B_S$.
  In view of the last argument of Section 1, there exist two schemes as:
  \[ a=a_1s_1,~ a_1t_1=a_2s_2,~ a_2t_2=a_3s_3, \ldots,~a_nt_n=a'\] and
  \[ b=b_1u_1,~ b_1v_1=b_2u_2,~ b_2v_2=b_3u_3, \ldots,~b_nv_n=b'\]
 both of length $n$ for some $n \in \mathbb{N},a_i\in A_S~,b_i \in B_S,~1\leq i \leq n$. By assumption, there exists a scheme connecting $(s_1,u_1)$ to $(t_1,v_1)$ in $S\times S$ of the form
\begin{multline*}
 (s_1,u_1)=(x_1,y_1)w_1,~ (x_1,y_1)z_1=(x_2,y_2)w_2,\\ (x_2,y_2)z_2=(x_3,y_3)w_3, \ldots,~(x_m,y_m)z_m=(t_1,v_1)
 \end{multline*}
  which yields the scheme
   \begin{multline*}
 (a,b)=(a_1s_1,b_1u_1)=(a_1x_1,b_1y_1)w_1,~ (a_1x_1,b_1y_1)z_1=(a_1x_2,b_1y_2)w_2,\\ (a_1x_2,b_1y_2)z_2=(a_1x_3,b_1y_3)w_3, \ldots,~(a_1x_m,b_1y_m)z_m=(a_1t_1,b_1v_1)
 \end{multline*} and hence we can assert that $(a,b)$ and $(a_1t_1, b_1v_1)$ are in the same indecomposable component. Processing inductively, we conclude that $(a,b)$ and $(a_nt_n,b_nv_n)=(a',b')$ are in the same indecomposable component.

 \epf

Regarding Theorem \ref{th1} and Corollary \ref{co4} the next corollary is obtained.
\bco \lb{co} For a finite monoid $S$, $S^{S\times S}$ is indecomposable if and only if $S\times S$ is indecomposable. \eco
If products of indecomposable acts are indecomposable, then $S^I$ is indecomposable for each nonempty set $I$, though, in comparison with Theorem \ref{th1}, this is a strict implication (see Example \ref{ex1}). Hereby, we need an additional condition on $S$ to fill the gap namely \emph{Condition right(left)-FI} under which there exists a fixed natural number $n$ such that any pair of elements in any indecomposable right(left) $S$-act can be connected via a scheme of length $n$ (see \cite[Corollary 2.11]{Bul}).

 In the next proposition we characterize monoids satisfying Condition right-FI.

\bpr \lb{pr4} Monoids satisfying condition right-FI are precisely left reversible monoids.\epr

\bpf
\textit{Necessity}. Suppose, by way of contradiction,  that $aS\cap bS=\emptyset$ for some $a,b\in S$. For each $i\in \mathbb{N}$, set $S_i=\{(i,s)~|~s\in S\}$ and endow it with the right $S$-action $(i,s)t=(i,st)$ for $s,t\in S$. Let us denote the element $(i,s)$ by $s^{(i)}$ for $i\in \mathbb{N}, s\in S$. For $n\in \mathbb{N}$ we define $A_n=\bigcup\limits_{i=1}^n S_i$ and  $\bar{A_n}= A_n/\rho_n$ where $\rho_n$ is the right congruence on $A_n$ generated by the pairs $(a^{(i)},b^{(i+1)}),~1\leq i \leq n$. Because of $aS \cap bS=\emptyset$, for $x,y\in A_n$ we have $x \rho_n y$ only if $x,y \in S_i \cup S_{i+1}$ for some $1\leq i\leq n-1$. On the other hand, since $a$ and $b$ are not right invertible, ${[1^{(i)}]}_{\rho_n}=\{1^{(i)}\}$ for any $1\leq i\leq n$. Now we have the following scheme
\begin{multline*} {[1^{(1)}]}_{\rho_n}={[1^{(1)}]}_{\rho_n}1,~{[1^{(1)}]}_{\rho_n}a={[1^{(2)}]}_{\rho_n}b,\\ {[1^{(2)}]}_{\rho_n}a={[1^{(3)}]}_{\rho_n}b
~, \ldots,~ {[1^{(n-1)}]}_{\rho_n}a={[1^{(n)}]}_{\rho_n}b \end{multline*} of length $n$ connecting ${[1^{(1)}]}_{\rho_n}$ to ${[1^{(n)}]}_{\rho_n}$. Since $\{{[1^{(i)}]}_{\rho_n}~|~1\leq i \leq n \}$ is a generating set for $\bar{A_n}$ and these generators are all in the same indecomposable component, $\bar{A_n}$ is indecomposable. Let there exist another scheme
 \begin{multline*} {[1^{(1)}]}_{\rho_n}={[a_1]}_{\rho_n}s_1,~ ~~{[a_1]}_{\rho_n}t_1={[a_2]}_{\rho_n}s_2, \\~ {[a_2]}_{\rho_n}t_2={[a_3]}_{\rho_n}s_3, \ldots,~{[a_m]}_{\rho_n}t_m={[1^{(n)}]}_{\rho_n} \end{multline*}
 of length $m$ connecting ${[1^{(1)}]}_{\rho_n}$ to ${[1^{(n)}]}_{\rho_n}$ for $s_i,t_i\in S,~a_i\in A_n,~1\leq i \leq m$. From  $1^{(1)} \rho_n ~ a_1s_1 $ and $a_mt_m ~\rho_n   ~1^{(n)}$, we get $a_1s_1=1^{(1)}$ and  $a_m t_m=1^{(n)}$ which imply that $a_1\in S_1$ and $a_m\in S_n$. Since $a_1t_1~\rho_n ~ a_2s_2$, $a_2\in S_1\cup S_2$. Continuing inductively, $a_m\in S_1 \cup S_2 \ldots \cup S_m$. Now, $a_m\in S_n$ implies that $n \leq m$ and hence, the shortest scheme connecting ${[1^{(1)}]}_{\rho_n}$ to ${[1^{(n)}]}_{\rho_n}$ is of length $n$. Considering $\bar{A_n}$,  for each $n\in \mathbb{N}$, $S$ doesn't satisfy condition right-FI, a contradiction.

 \textit{Sufficiency}. If S is left reversible then, by Proposition \ref{pr1}, any pair of elements in any indecomposable right $S$-act is connected by a scheme of length 2.

\epf

In the next proposition we characterize monoids for which products of indecomposable acts are indecomposable.

 \bpr \lb{pr5} For a monoid $S$ the following are equivalent:

$~~i)$ products of indecomposable right $S$-acts are indecomposable,

$~ii)$ $S$ is left reversible and $S^{S\times S}$ is indecomposable,

$iii)$ $S$ satisfies condition right-FI  and $S^{S\times S}$ is indecomposable,

$\,iv)$ non-zero cofree acts are decomposable and $S^{S\times S}$ is indecomposable,

$\,\,\,v)$ all subacts of indecomposable right $S$-acts are indecomposable and $S^{S\times S}$ is indecomposable.
\epr

\bpf By virtue of Propositions \ref{pr2}, \ref{pr3} and \ref{pr4} it is enough to prove that the first two statements are equivalent.

$i\Longrightarrow ii$ Suppose, contrary to our claim, that $S$ is not left reversible. For each $n\in \mathbb{N}$, let $\bar{A_n}$ be the right $S$-act constructed in the proof of Proposition \ref{pr4}. Set $A=\prod \limits _{n\in \mathbb{N}}\bar{A_n}$ which is indecomposable by assumption. Therefore there is a scheme of length $m$, connecting ${({[1^{(1)}]}_{\rho_n})}_{n\in \mathbb{N}}$ to ${({[1^{(n)}]}_{\rho_n})}_{n\in \mathbb{N}}$. Considering this scheme componentwise, for each $n\in \mathbb{N}$ there exists a scheme of length $m$ in $\bar{A_n}$, connecting $ {[1^{(1)}]}_{\rho_n}$ to ${[1^{(n)}]}_{\rho_n}$. But according to the proof of Proposition \ref{pr4}, $m\geq n$ for each $n\in \mathbb{N}$, a contradiction.

$ii)\Longrightarrow i)$ Let $S$ be a left reversible monoid and let $\{A_i~|~i\in I \}$ be a family of indecomposable right $S$-acts and ${(a_i)}_{i\in I},~{(b_i)}_{i\in I}\in \prod \limits _{i\in I} A_i$. In light of Proposition \ref{pr1}, for each $i\in I$, there exist $s_i,t_i \in S$ such that $a_is_i=b_it_i$. Let us denote a typical element of $S^I$, with the same element $x\in S$ in each component, by $(x)_{i\in I}$. Since $S^I$ is indecomposable by Corollary \ref{co4} and ${(s_i)}_{i\in I},~{(1)}_{i\in I}\in S^I$, another application of Proposition \ref{pr1}, part $ii$, provides $s,t \in S$ such that ${(s_i)}_{i\in I}s={(1)}_{i\in I}t$. The same arguments provides existence of $\alpha, \beta \in S$ such that ${(t_is)}_{i\in I}\alpha=~{(t)}_{i\in I}\beta$. Now for each $i\in I$ we have $ a_it\alpha =a_is_is\alpha=b_it_is\alpha=b_it\beta$, which yields
 ${(a_i)_{i\in I}}{t\alpha}={(b_i)_{i\in I}}{t\beta}$ as desired.
\epf
For commutative monoids, the left reversibility condition in Proposition \ref{pr5} is fulfilled and the following corollary is obtained.

\bco \lb{co5} For a commutative monoid $S$ products of indecomposable acts are indecomposable if and only if $S^{S\times S}$ is indecomposable. \eco

Recall that a monoid $S$ is called right collapsible if for any $s,t\in S$ there exists $u\in S$ such that $su=tu$.

\ble \lb{le6} For a left reversible monoid $S$, finite products of indecomposable right $S$-acts are indecomposable if and only if $S$ is right collapsible.
\ele
\bpf Let $S$ be a left reversible monoid, $s,t \in S$, and let finite products of indecomposable right $S$-acts be indecomposable. Since $S\times S$ is indecomposable, using Proposition \ref{pr1} there exist $u,v \in S$ such that $(1,s)u=(1,t)v$ that is $u=v$ and $su=tu$.

Conversely, suppose that $S$ is a right collapsible monoid and $(a,b),(c,d) \in S \times S$. Under our assumption, there exist $u_1,u_2 \in S$ such that $au_1=cu_1,~bu_2=du_2$. Also $u_1u=u_2u$ for some $u \in S$. Then  $(a,b)u_1u=(c,d)u_1u$ as desired.

\epf
Considering the strict implication $right~collapsible \Longrightarrow left~reversible$ for monoids, Lemma \ref{le6} and Theorem \ref{th1} give the following result.

\bco \lb{co6} A monoid $S$ is right collapsible if and only if $S$ is left reversible and $S\times S$ is indecomposable.
\eco

\bte \lb{th2}For a monoid $S$ products of indecomposable right $S$-acts are indecomposable if and only if $S$ has a right zero.
\ete
\bpf \textit{Necessity}. Suppose that $S$ is represented by an index set $I$ as $S=\{s_i~|~i\in I\}$. In view of Proposition \ref{pr5} and Corollary \ref{co4}, $S$ is left reversible and $S^I$ is indecomposable. By Proposition \ref{pr1}, ${(s_i)}_{i\in I}s= {(1)}_{i\in I}t$ for some $s,t\in S$. Thus for each $x\in S$, $xs=t$. Taking $x=1$, gives $s=t$ that is $xt=t$ for every $ x\in S$ as desired.

\textit{Sufficiency}. Let $z\in S$ be a right zero and let $I$ be a nonempty set. Take ${(a_i)}_{i\in I},{(b_i)}_{i\in I} \in S^I$. Thus ${(a_i)}_{i\in I}z={(b_i)}_{i\in I}z$ and then $S^I$ is indecomposable. By assumption $S$ is left reversible and regarding Proposition \ref{pr5} the result follows.

\epf

\brm It is worth to mention that the two conditions in the sufficient part of Corollary \ref{co6}, regarding Example \ref{ex1}, are independent. Besides, analogously to the Corollary \ref{co6} and the strict implication {\it monoid with right zero $\Longrightarrow$ left reversible monoid}, Theorems \ref{th1} and \ref{th2} state that monoids with right zeros are precisely left reversible monoids for which $S^{S\times S}$ is indecomposable.

\erm

The next example shows that the conditions occurring in part $ii)$ of Proposition \ref{pr5} are independent.

\bex \lb{ex1}
Let $S$ be a left reversible monoid which doesn't have a right zero for instance a nontrivial  finite group. Proposition \ref{pr5} and \ref{th2} imply that $S^{S\times S}$ is not indecomposable. On the other hand, let $S=\mathcal{T}_n$ consist of transformations of the set $\{1,2,\ldots, n\}$ with mappings written on the left side for some $1\neq n \in \mathbb{N}$. Set $c_i=\left ( \begin{array}{cccc} 1&2&\ldots&n\\
i&i&\ldots&i \end{array}\right )$ for $1\leq i\leq n$. For $(\alpha,\beta) \in S \times S, $ we have $(\alpha,\beta)c_1= (c_{i_1},c_{i_2})$ where $i_1,i_2\in \{1,2,\ldots, n\}$  which states that $(\alpha,\beta)$ and $(c_{i_1},c_{i_2})$ are in the same indecomposable component. Now the following scheme
\begin{alignat*}{3}
~~~~~~~~~~~~~~~~~~~(c_1,c_1)&=(1,c_1)c_1 &&~~,~~&(1,c_1)c_{i_1} &=(c_{i_1},c_1)1\\
(c_{i_1},c_1)1&=(c_{i_1},1)c_1 &&~~,~~&(c_{i_1},1)c_{i_2}&=(c_{i_1},c_{i_2})
\end{alignat*}

in $S\times S$, implies that $(c_1,c_1)$ and $(c_{i_1},c_{i_2})$ are in the same indecomposable component. From this $S\times S$ is indecomposable. Using Corollary \ref{co} $S^I$ is indecomposable for each nonempty set $I$ and since $S$ contains two left zeros, $S$ is not left reversible.
\eex
In the above example we observed that for the monoid $S=\mathcal{T}_n$, $S\times S$ is indecomposable. So a question can be posed that \\{\bf whether for the monoid of full transformations of a nonempty set $X$, $S^I$ is indecomposable for each nonempty set $I$}.

 Note that in \cite{Bul}, Proposition 3.8 states that for a proper right ideal $K$ of a monoid $S$ if the Rees factor act $S/K$ is finitely product flat then $S/K$ is super flat. So a natural question that comes to the mind is the case that $K=S$. In the next proposition we show that in this case product flatness is equivalent to super flatness. Indeed  in \cite[Corollary 2.11]{Bul} it is proved that the one element right $S$-act $\Theta_S$ is product flat if and only if $S$ satisfies condition left-$FI$ and the left $S$-act $_SS^I$ is indecomposable for each nonempty set $I$. Hereby we give the next proposition which is an improvement of this result.
\bpr \lb{pr8} For a monoid $S$ the following are equivalent:

$~~i)$ the one element right $S$-act $\Theta_S$ is super flat,

$~ii)$ the one element right $S$-act $\Theta_S$ is product flat,

$iii)$ products of indecomposable left $S$-acts are indecomposable.

$\,iv)$ $S$ contains a left zero.
\epr
\bpf $i \Longrightarrow ii$ is trivial. $ii \Longrightarrow iii$ follows immediately by Proposition \ref{pr5} and \cite[Corollary 2.11]{Bul}. Theorem \ref{th2} implies the equivalence of $iii$ and $iv$  and $iv \Longrightarrow i$ follows by \cite[Corollary 3.6]{Bul}.
\epf
\begin{center}\section*{Acknowledgements} \end{center}
The authors wish to express their appreciation to the anonymous referee for his/her constructive input into the paper.


\begin{thebibliography}{99}
\bibitem{adamek} Adamek, J.,  Herrlich, H. and Strecker, G. \textit{ Abstract and Concrete Categories The Joy of Cats} (John Wiley and Sons, New York, 1990).
\bibitem{Bul1} Bulman-Fleming, S. \textit{Products
of projective S-systems}, Comm. Algebra \textbf{19} (3), 951--964, 1991.
\bibitem{Bul} Bulman-Fleming, S. and Laan, V. \textit{Tensor products and preservation of limits, for acts over monoids}, Semigroup Forum \textbf{63}, 161--179, 2001.
\bibitem{Bul2} Bulman-Fleming, S and McDowell, K.   \textit{Coherent monoids}, in: Latices, Semigroups and Universal Algebra,
ed. J. Almeida et al. (Plenum Press, New York, 1990).
\bibitem{Gould} Gould, V. \textit{Coherent monoids}, J. Austral. Math. Soc. \textbf{53}(Series A), 166--182, 1992.
\bibitem{Kilp} Kilp, M., Knauer, U. and Mikhalev, A. \textit{Monoids, Acts and Categories}, (W. de gruyter, Berlin, 2000).
\bibitem{niko} Nico, W.R.  \textit{A classification of indecomposable $S$-sets}, J. Algebra \textbf{54}(1), 260--272, 1978.
\bibitem{Ren} Renshaw, J. \textit{Monoids for which condition (P)acts are projective}, Semigroup Forum \textbf{61}(1), 46--56, 1998.
\bibitem{sedaghat} Sedaghatjoo, M., Khosravi, R. and Ershad, M. \textit{Principally weakly and weakly coherent monoids}, Comm. Algebra \textbf{37}(12), 4281--4295, 2009.


\end{thebibliography}
\end{document}